%% file: main.tex
\documentclass{iise}

\usepackage{bbold}
\usepackage{amsfonts,amsmath,amssymb,amsthm}
\usepackage{pifont}
\usepackage{accents}
\usepackage{mathrsfs}
\usepackage[dvipsnames]{xcolor}
\usepackage{pdflscape}
\usepackage{mathtools}
\usepackage{cancel}
\usepackage{hyperref}
\usepackage{nicefrac}

\usepackage{booktabs}
\usepackage{todonotes}

\usepackage{multirow}

\usepackage{float}
\floatstyle{ruled}
\newfloat{model}{thp}{lop}
\floatname{model}{Model}

\theoremstyle{definition}

\usepackage{algorithmic}
\usepackage[ruled,linesnumbered,noresetcount,vlined]{algorithm2e}

\newcommand{\nsfai}{{\fontsize{11}{12}\selectfont\normalfont NSF Artificial Intelligence Institute for Advances in Optimization}}

\newcommand{\gatech}{{\fontsize{11}{12}\selectfont\normalfont Georgia Institute of Technology, Atlanta, GA, USA}}

\title{\titlesize Dual Conic Proxy for Semidefinite Relaxation of\\ AC Optimal Power Flow}

\author{Guancheng Qiu, Mathieu Tanneau, Pascal Van Hentenryck \\
\nsfai \\ \gatech}
\authorlist{Qiu, Tanneau and Van Hentenryck}

\begin{document}

\pagestyle{headings}

\maketitle

\begin{abstract}
{\small
The nonlinear, non-convex AC Optimal Power Flow (AC-OPF) problem is fundamental for power systems operations.
The intrinsic complexity of AC-OPF has fueled a growing interest in the development of optimization proxies for the problem, i.e., machine learning models that predict high-quality, close-to-optimal solutions.
More recently, dual conic proxy architectures have been proposed, which combine machine learning and convex relaxations of AC-OPF, to provide valid certificates of optimality using learning-based methods.
Building on this methodology, this paper proposes, for the first time, a dual conic proxy architecture for the semidefinite (SDP) relaxation of AC-OPF problems.
Although the SDP relaxation is stronger than the second-order cone relaxation considered in previous work, its practical use has been hindered by its computational cost.
The proposed method combines a neural network with a differentiable dual completion strategy that leverages the structure of the dual SDP problem.
This approach guarantees dual feasibility, and therefore valid dual bounds, while providing orders of magnitude of speedups compared to interior-point algorithms.
The paper also leverages self-supervised learning, which alleviates the need for time-consuming data generation and allows to train the proposed models efficiently.
Numerical experiments are presented on several power grid benchmarks with up to 500 buses.
The results demonstrate that the proposed SDP-based proxies can outperform weaker conic relaxations, while providing several orders of magnitude speedups compared to a state-of-the-art interior-point SDP solver.
}
\end{abstract}

\section*{Keywords}
AC optimal power flow, convex relaxation, semidefinite programming, neural network, self-supervised learning

\input{tex/definitions}

\input{tex/introduction}
\input{tex/formulations}
\input{tex/learning}

\input{tex/results}
\input{tex/conclusion}

\section*{Acknowledgments}
This research was supported by NSF award 2112533 and ARPA-E PERFORM award AR0001136.

\bibliography{ref}

\end{document}

%% file: tex/definitions.tex
\makeatletter
\DeclareRobustCommand{\cev}[1]{%
  {\mathpalette\do@cev{#1}}%
}
\newcommand{\do@cev}[2]{%
  \vbox{\offinterlineskip
    \sbox\z@{$\m@th#1 x$}%
    \ialign{##\cr
      \hidewidth\reflectbox{$\m@th#1\vec{}\mkern4mu$}\hidewidth\cr
      \noalign{\kern-\ht\z@}
      $\m@th#1#2$\cr
    }%
  }%
}

\newcommand{\ubar}[1]{\underaccent{\bar}{#1}}

\makeatother

\newcommand{\clarabelQuad}{Clarabel\textsubscript{quad}}

\newcommand{\im}{\mathbf{j}}

\newcommand{\trace}{\operatorname{tr}}

\newcommand{\gff}{g^{\text{ff}}}
\newcommand{\gft}{g^{\text{ft}}}
\newcommand{\gtf}{g^{\text{tf}}}
\newcommand{\gtt}{g^{\text{tt}}}
\newcommand{\bff}{b^{\text{ff}}}
\newcommand{\bft}{b^{\text{ft}}}
\newcommand{\btf}{b^{\text{tf}}}
\newcommand{\btt}{b^{\text{tt}}}
\newcommand{\yff}{Y^{\text{ff}}}
\newcommand{\yft}{Y^{\text{ft}}}
\newcommand{\ytf}{Y^{\text{tf}}}
\newcommand{\ytt}{Y^{\text{tt}}}

\newcommand{\pg}{\mathbf{p}^{\text{g}}}
\newcommand{\qg}{\mathbf{q}^{\text{g}}}
\newcommand{\Sg}{\mathbf{S}^{\text{g}}}

\newcommand{\vm}{\mathbf{v}}
\newcommand{\va}{\mathbf{\theta}}
\newcommand{\V}{\mathbf{V}}
\newcommand{\vr}{\mathbf{e}}
\newcommand{\vi}{\mathbf{f}}

\newcommand{\pf}{\vec{\mathbf{p}}}
\newcommand{\qf}{\vec{\mathbf{q}}}
\newcommand{\Sf}{\vec{\mathbf{S}}}
\newcommand{\pt}{\cev{\mathbf{p}}}
\newcommand{\qt}{\cev{\mathbf{q}}}
\newcommand{\St}{\cev{\mathbf{S}}}

\newcommand{\wmsoc}{\mathbf{w}}
\newcommand{\wcsoc}{\mathbf{w}^{\text{re}}}
\newcommand{\wssoc}{\mathbf{w}^{\text{im}}}

\newcommand{\Wsdp}{\mathbf{W}}
\newcommand{\Xsdp}{\mathbf{X}}
\newcommand{\Wrsdp}{\mathbf{W}^{\text{re}}}
\newcommand{\Wisdp}{\mathbf{W}^{\text{im}}}

\DeclarePairedDelimiterX{\inner}[2]{\langle}{\rangle}{#1, #2}
\newcommand{\symmat}{\mathbb{S}}
\newcommand{\skewmat}{\overline{\symmat}}
\newcommand{\psdsymmat}{\symmat_{+}}
\newcommand{\hermat}{\mathbb{H}}
\newcommand{\psdhermat}{\hermat_{+}}
\newcommand{\sympart}{\mathcal{H}}  
\newcommand{\skewpart}{\mathcal{Z}}  
\newcommand{\blocksymmat}{\mathbb{S}^{B}}
\newcommand{\blockpsdsymmat}{\mathbb{S}^{B}_{+}}
\newcommand{\E}{\mathbf{E}}
\newcommand{\Ere}{\E^{\text{re}}}
\newcommand{\Eim}{\E^{\text{im}}}

\newcommand{\pd}{\text{p}^{\text{d}}}
\newcommand{\qd}{\text{q}^{\text{d}}}
\newcommand{\Sd}{\text{S}^{\text{d}}}
\newcommand{\Sgmin}{\underline{\text{S}}^{\text{g}}}
\newcommand{\Sgmax}{\overline{\text{S}}^{\text{g}}}
\newcommand{\pgmin}{\underline{\text{p}}^{\text{g}}}
\newcommand{\pgmax}{\overline{\text{p}}^{\text{g}}}
\newcommand{\qgmin}{\underline{\text{q}}^{\text{g}}}
\newcommand{\qgmax}{\overline{\text{q}}^{\text{g}}}
\newcommand{\vmmin}{\underline{\text{v}}}
\newcommand{\vmmax}{\overline{\text{v}}}
\newcommand{\dvamin}{\underline{\Delta \va}_{ij}}
\newcommand{\dvamax}{\overline{\Delta \va}_{ij}}
\newcommand{\wcsocmin}{\underline{\mathbf{w}}^{\text{re}}}
\newcommand{\wcsocmax}{\overline{\mathbf{w}}^{\text{re}}}
\newcommand{\wssocmin}{\underline{\mathbf{w}}^{\text{im}}}
\newcommand{\wssocmax}{\overline{\mathbf{w}}^{\text{im}}}

\newcommand{\gammaPwfr}{\vec{\gamma}^{\text{p,w}}}
\newcommand{\gammaQwfr}{\vec{\gamma}^{\text{q,w}}}
\newcommand{\gammaPwto}{\cev{\gamma}^{\text{p,w}}}
\newcommand{\gammaQwto}{\cev{\gamma}^{\text{q,w}}}
\newcommand{\gammaPrfr}{\vec{\gamma}^{\text{p,r}}}
\newcommand{\gammaQrfr}{\vec{\gamma}^{\text{q,r}}}
\newcommand{\gammaPrto}{\cev{\gamma}^{\text{p,r}}}
\newcommand{\gammaQrto}{\cev{\gamma}^{\text{q,r}}}
\newcommand{\gammaPifr}{\vec{\gamma}^{\text{p,i}}}
\newcommand{\gammaQifr}{\vec{\gamma}^{\text{q,i}}}
\newcommand{\gammaPito}{\cev{\gamma}^{\text{p,i}}}
\newcommand{\gammaQito}{\cev{\gamma}^{\text{q,i}}}

\newcommand{\Y}{{\mathbf{Y}}}
\newcommand{\G}{\mathbf{G}}
\newcommand{\B}{\mathbf{B}}

\newcommand{\M}{\mathbf{M}}
\newcommand{\Mfr}{\vec{\M}}
\newcommand{\Mto}{\cev{\M}}

\newcommand{\lambdaP}{\lambda^{\text{p}}}
\newcommand{\lambdaQ}{\lambda^{\text{q}}}
\newcommand{\lambdaPf}{\vec{\phantom{\lambda}}\hspace{-0.58em}{\lambda}^{\text{p}}}
\newcommand{\lambdaPt}{\cev{\phantom{\lambda}}\hspace{-0.35em}{\lambda}^{\text{p}}}
\newcommand{\lambdaQf}{\vec{\phantom{\lambda}}\hspace{-0.58em}{\lambda}^{\text{q}}}
\newcommand{\lambdaQt}{\cev{\phantom{\lambda}}\hspace{-0.35em}{\lambda}^{\text{q}}}
\newcommand{\nuThermalfr}{\vec{\nu}}
\newcommand{\nuThermalSfr}{\vec{\nu}^{\, \text{s}}}
\newcommand{\nuThermalPfr}{\vec{\nu}^{\, \text{p}}}
\newcommand{\nuThermalQfr}{\vec{\nu}^{\, \text{q}}}
\newcommand{\nuThermalto}{\cev{\nu}}
\newcommand{\nuThermalSto}{\cev{\nu}^{\, \text{s}}}
\newcommand{\nuThermalPto}{\cev{\nu}^{\, \text{p}}}
\newcommand{\nuThermalQto}{\cev{\nu}^{\, \text{q}}}
\newcommand{\nuWprod}{\nu^{\text{w}}}
\newcommand{\omegaf}{\omega^{\text{f}}}
\newcommand{\omegat}{\omega^{\text{t}}}
\newcommand{\omegar}{\omega^{\text{re}}}
\newcommand{\omegai}{\omega^{\text{im}}}
\newcommand{\muPg}{{\mu}^{\text{pg}}}
\newcommand{\muPgMin}{\ubar{\mu}^{\text{pg}}}
\newcommand{\muPgMax}{\bar{\mu}^{\text{pg}}}
\newcommand{\muQg}{{\mu}^{\text{qg}}}
\newcommand{\muQgMin}{\ubar{\mu}^{\text{qg}}}
\newcommand{\muQgMax}{\bar{\mu}^{\text{qg}}}
\newcommand{\muWm}{{\mu}^{\text{w}}}
\newcommand{\muWmMin}{\ubar{\mu}^{\text{w}}}
\newcommand{\muWmMax}{\bar{\mu}^{\text{w}}}
\newcommand{\muAngleDiff}{{\mu}^{\Delta\theta}}
\newcommand{\muAngleDiffMin}{\ubar{\mu}^{\Delta\theta}}
\newcommand{\muAngleDiffMax}{\bar{\mu}^{\Delta\theta}}
\newcommand{\muWr}{{\mu}^{\text{wr}}}
\newcommand{\muWrMin}{\ubar{\mu}^{\text{wr}}}
\newcommand{\muWrMax}{\bar{\mu}^{\text{wr}}}
\newcommand{\muWi}{{\mu}^{\text{wi}}}
\newcommand{\muWiMin}{\ubar{\mu}^{\text{wi}}}
\newcommand{\muWiMax}{\bar{\mu}^{\text{wi}}}

\newcommand{\Ssdp}{\mathbf{S}}

\newcommand{\pdRef}{\bar{\text{p}}^{\text{d}}}
\newcommand{\qdRef}{\bar{\text{q}}^{\text{d}}}

%% file: tex/introduction.tex
\section{Introduction}
\label{sec:introduction}
The AC Optimal Power Flow (AC-OPF) problem aims to find the most cost-effective to dispatch power generation while satisfying demand and engineering constraints.
Despite its importance to power system operation, planning and electricity markets, its complex, nonlinear, and non-convex nature makes it challenging to solve and limits its practical use.

The need to repeatedly solve similar instances has led to the development of optimization proxies, i.e., machine learning (ML) models which can approximate the input-output mapping of AC-OPF solvers in milliseconds.
Most existing proxies for DC- and AC-OPF focus on predicting primal solutions without providing optimality guarantees \cite{fioretto2020predicting,donti2021dc3,huang2021deepopf,Chen2023_E2ELR}.
Hence, these approaches may provide sub-optimal solutions, which represents significant economic losses, and are therefore not sufficient for electricity market-clearing applications, which require near-optimal solutions \cite{MISO_BPM_002}.

To address this limitation, Dual Conic Proxies (DCP) were introduced in \cite{qiu2024dual}, which proposes to learn to generate valid lower bounds by predicting feasible dual solutions for a second-order cone (SOC) relaxation of AC-OPF proposed by Jabr \cite{Jabr2006_SOCRelaxationOPF}.
A more general Dual Lagrangian Learning framework was then proposed in   \cite{tanneau2024dual} for conic optimization problems, with results reported on linear and SOC problems.
While the approach of \cite{qiu2024dual} scales to large systems, it is intrinsically limited by the strength of the Jabr relaxation.
The reader is referred to \cite{Molzahn2019_OPF_survey} and to \cite{qiu2024dual,tanneau2024dual} for a detailed review of convex relaxations of AC-OPF, and of dual optimization proxies, respectively.

It is well-known that semidefinite (SDP) formulations provide strong relaxations for OPF problems \cite{kocuk2016strong,Molzahn2019_OPF_survey}, at the price of significant computational and numerical challenges.
As a first step towards achieving practical DCP methodology for SDP problems, this work develops, for the first time, a 
DCP methodology for an SDP relaxation of AC-OPF.
The paper presents a new dual-feasible architecture for this SDP relaxation, and conducts numerical results on small-to-medium power grids that demonstrate the benefits of the proposed approach.

%% file: tex/formulations.tex
\section{AC-OPF and a Semidefinite Relaxation}
\label{sec:formulation}

The imaginary unit is denoted by $\im$, i.e., $\im^{2} \, {=} \, -1$.
The complex conjugate of $z \, {\in} \, \mathbb{C}$ is $z^{\star}$.
$\Re(\cdot)$ and $\Im(\cdot)$ denotes the real and imaginary parts of a complex number.
In all that follows, let $\E_{ij}$ denote a square matrix of appropriate dimension, whose $(i, j)$ entry is equal to $1$, and all other zeros.
The identity matrix is denoted by $I$.
The symmetric and skew-symmetric parts of square matrix $\mathbf{A}$ are denoted by $\mathbf{A}^{+} = (\mathbf{A} + \mathbf{A}^{\top})/2$ and $\mathbf{A}^{-} = (\mathbf{A} - \mathbf{A}^{\top})/2$.
Note that $\mathbf{A} = \mathbf{A}^{+} + \mathbf{A}^{-}$.
The smallest eigenvalue of a Hermitian matrix $\mathbf{X}$ is denoted by $\lambda_{min}(\mathbf{X})$.
The cones of positive semidefinite Hermitian (resp. symmetric) matrices of order $n$ is denoted by $\mathbb{H}_{+}^{n}$ (resp. $\mathbb{S}^{n}_{+}$),
and the second-order cone of order $n$ is denoted by $\mathcal{Q}^{n} = \{x \in \mathbb{R}^{n} \, | \, x_{1} \geq \sqrt{x_{2}^{2} + ... + x_{n}^{2}} \}$.

The paper considers a power grid, represented as a simple directed graph, whose sets of buses and branches are denoted by $\mathcal{N}$ and $\mathcal{E}$, respectively.
Each branch is represented as a directed edge $(i, j)$ from bus $i$ to $j$, with admittance matrix
\begin{align}
    Y_{ij} =
    \begin{pmatrix}
        \yff_{ij} & \yft_{ij}\\
        \ytf_{ij} & \ytt_{ij}
    \end{pmatrix}
    = 
    \begin{pmatrix}
        \gff_{ij} + \im \bff_{ij} & \gft_{ij} + \im \bft_{ij}\\
        \gtf_{ij} + \im \btf_{ij} & \gtt_{ij} + \im \btt_{ij}
    \end{pmatrix}
    \in 
    \mathbb{C}^{2 \times 2}.
\end{align}
The shunt admittance and power demand at node $i$ are denoted by $Y^{s}_{i} = g^{s}_{i} + \im b^{s}_{i} \in \mathbb{C}$ and $\Sd_{i} = \pd_{i} + \im \qd_{i} \in \mathbb{C}$, respectively.
For ease of reading and without loss of generality, we formulate the problems assuming that exactly one generator is connected to each bus, and that generation costs are linear.

\subsection{The AC-Optimal Power Flow Formulation}

The formulation of AC-OPF considered in this paper is detailed in Model~\ref{model:AC-OPF}.
The decision variables comprise complex generation $\Sg = \pg + \im \qg$,
complex nodal voltage $\V = \vm \angle \va$,
and forward and reverse power flows $\Sf = \pf + \im \qf$ and $\St = \pt + \im \qt$.
The objective \eqref{eq:ACOPF:objective} minimizes total generation costs.
Constraint \eqref{eq:ACOPF:kirchhoff} enforces Kirchhoff current law at each bus.
Constraints \eqref{eq:ACOPF:ohm_fr}--\eqref{eq:ACOPF:ohm_to} express power flows on each branch using Ohm's law,
and constraints \eqref{eq:ACOPF:thermal_limits} enforce thermal limits on forward and reverse power flows.
Finally, constraints \eqref{eq:ACOPF:voltage_bounds} and \eqref{eq:ACOPF:reactive_dispatch_bounds} enforce minimum and maximum limits on nodal voltage magnitude and power generation.

\begin{figure}[!t]
    \centering
    \begin{minipage}[t]{0.47\textwidth}
        \begin{model}[H]
            \caption{The AC-OPF model}
            \label{model:AC-OPF}
            \begin{subequations}
            \footnotesize
            \label{eq:ACOPF}
            \begin{align}
                \min \quad 
                & \sum_{i \in \mathcal{N}} c_{i} \pg_{i}
                    \label{eq:ACOPF:objective}
                    \\
                \textrm{s.t.} \quad
                & \Sg_{i} - \Sd_{i} - Y_{i}^{s^{*}} |\V_{i}|^{2} = \sum_{(i,j) \in \mathcal{E}} \Sf_{ij} + \sum_{(j,i) \in \mathcal{E}} \St_{ji}
                    && \forall i \in \mathcal{N}
                    \label{eq:ACOPF:kirchhoff}
                    \\
                & \Sf_{ij} = {\yff_{ij}}^{*} |\V_{i}|^{2} - {\yft_{ij}}^{*} \V_{i} \V_{j}^{*} 
                    && \forall ij \in \mathcal{E}
                    \label{eq:ACOPF:ohm_fr}
                    \\
                & \St_{ij} = {\ytt_{ij}}^{*} |\V_{j}|^{2} - {\ytf_{ij}}^{*} \V_{i}^{*} \V_{j}
                    && \forall ij \in \mathcal{E}
                    \label{eq:ACOPF:ohm_to}
                    \\
                & |\Sf_{ij}|, |\St_{ij}| \leq \bar{s}_{ij}
                    && \forall ij \in \mathcal{E}
                    \label{eq:ACOPF:thermal_limits}
                    \\
                & \vmmin_{i} \leq |\V_{i}| \leq \vmmax_{i} 
                    && \forall i \in \mathcal{N}
                    \label{eq:ACOPF:voltage_bounds}
                    \\
                & \Sgmin_{i} \leq \Sg_{i} \leq \Sgmax_{i}
                    && \forall i \in \mathcal{N}
                    \label{eq:ACOPF:reactive_dispatch_bounds}
            \end{align}
            \end{subequations}
            \vspace{0.5em}
        \end{model}
    \end{minipage}
    \hfill
    \begin{minipage}[t]{0.45\textwidth}
        \begin{model}[H]
            \caption{The SDP-OPF model}
            \label{model:SDP-OPF}
            \begin{subequations}
            \footnotesize
            \begin{align}
                \min \quad 
                & \sum_{i \in \mathcal{G}} c_{i} \pg_{i}
                    \label{eq:SDPOPFC:obj}
                    \\
                \textrm{s.t.} \quad
                & \Sg_{i} - \Sd_{i} - Y_{i}^{s^{*}} \Wsdp_{ii} = \sum_{(i,j) \in \mathcal{E}} \Sf_{ij} + \sum_{(j,i) \in \mathcal{E}} \St_{ji}
                    && \forall i \in \mathcal{N}
                    \label{eq:SDPOPFC:kirchhoff}
                    \\
                & \Sf_{ij} = \inner{\yff_{ij} \E_{ii} + \yft_{ij} \E_{ij}}{\Wsdp}
                    && \forall ij \in \mathcal{E}
                    \\
                & \St_{ij} = \inner{\ytt_{ij} \E_{jj} + \ytf_{ij} \E_{ij}}{\Wsdp}
                    && \forall ij \in \mathcal{E}
                    \\
                & 
                    (\bar{s}_{ij}, \pf_{ij}, \qf_{ij}), (\bar{s}_{ij}, \pt_{ij}, \qt_{ij}) \in \mathcal{Q}^{3}
                    && \forall ij \in \mathcal{E}
                    \\
                & \vmmin_{i}^{2} \leq \Wsdp_{ii} \leq \vmmax_{i}^{2} 
                    && \forall i \in \mathcal{N} \\
                & \Sgmin_{i} \leq \Sg_{i} \leq \Sgmax_{i}
                    && \forall i \in \mathcal{N}
                    \label{eq:SDPOPFC:qg:bounds}
                    \\
                & \Wsdp \in \mathbb{H}^{|\mathcal{N}|}_{+}
                    \label{eq:SDPOPFC:psd}
            \end{align}
            \end{subequations}
        \end{model}
    \end{minipage}
\end{figure}

\subsection{A Semidefinite Relaxation of AC-OPF}

Introducing the change of variable $\Wsdp = \V \V^{\star}$ yields an equivalent formulation of AC-OPF if $\Wsdp \succeq 0$ and $\text{rank}(\Wsdp) = 1$.
Relaxing the rank constraint on $\Wsdp$ yields a semidefinite relaxation SDP-OPF, originally proposed in \cite{Bai2008_SDPRelaxationOPF}, which is stated in Model \ref{model:SDP-OPF}.
Because SDP-OPF is a relaxation of AC-OPF, its optimal value is a valid lower bound on the optimal value of AC-OPF.

Model \ref{model:DSDP-OPF} presents the conic dual of SDP-OPF, where $\mathcal{A}_{R}(\lambda, \mu, \nu)$ and $\mathcal{A}_{I}(\lambda, \mu, \nu)$ are defined as
\begin{align}
    \label{eq:AR_def}
    \mathcal{A}_{R}(\lambda, \mu, \nu) = 
        & \sum_{i \in \mathcal{N}} 
            (- g_{i}^{s} \lambdaP_{i} + b_{i}^{s} \lambdaQ_{i}) \E_{ii} 
            + (\muWmMin_{i} - \muWmMax_{i}) \E_{ii}\\
        &  + \sum_{(i, j) \in \mathcal{E}} \left(
                  \lambdaPf_{ij} \big( \gff_{ij} \E_{ii} + \gft_{ij} \E^{+}_{ij} \big)
                + \lambdaPt_{ij} \big( \gtt_{ij} \E_{jj} + \gtf_{ij} \E^{+}_{ij} \big)
                - \lambdaQf_{ij} \big( \bff_{ij} \E_{ii} + \bft_{ij} \E^{+}_{ij} \big)
                - \lambdaQt_{ij} \big( \btt_{ij} \E_{jj} + \btf_{ij} \E^{+}_{ij} \big)
                \right) \nonumber\\
    \mathcal{A}_{I}(\lambda, \mu, \nu) = 
        & \sum_{(i, j) \in \mathcal{E}} \left(
            \lambdaPf_{ij} \bft_{ij} \E^{-}_{ij}
            - \lambdaPt_{ij} \btf_{ij} \E^{-}_{ij}
            + \lambdaQf_{ij} \gft_{ij} \E^{-}_{ij}
            - \lambdaQt_{ij} \gtf_{ij} \E^{-}_{ij}
        \right).
\end{align}
By weak duality, any dual solution that satisfies \eqref{eq:DSDPOPF:pg}-\eqref{eq:DSDPOPF:psd} yields a valid lower bound on the optimal value of SDP-OPF and, in turn, the optimal value of AC-OPF.

\begin{model}[!t]
\caption{The DSDP-OPF model}
\label{model:DSDP-OPF}
\begin{subequations}
    \small
    \begin{align}
        \max_{\lambda, \mu, \nu} \quad 
        & \sum_{i \in \mathcal{N}} \left(
            \pd_{i} \lambdaP_{i}
            + \qd_{i} \lambdaQ_{i}
            + \vmmin_{i}^{2} \muWmMin_{i} - \vmmax_{i}^{2} \muWmMax_{i}
            + \pgmin_{i} \muPgMin_{i} - \pgmax_{i} \muPgMax_{i}
            + \qgmin_{i} \muQgMin_{i} - \qgmax_{i} \muQgMax_{i}
        \right)
        - \sum_{e \in \mathcal{E}} \bar{s}_{e} \left( \nuThermalSfr_{e} + \nuThermalSto_{e} \right)
        \label{eq:DSDPOPF:obj}
        \\
        \text{s.t.} \quad
        & 
            \lambdaP_{i} + \muPgMin_{i} - \muPgMax_{i} = c_{i}
            && \forall i \in \mathcal{N}
            \label{eq:DSDPOPF:pg}\\
        & 
            \lambdaQ_{i} + \muQgMin_{i} - \muQgMax_{i} = 0
            && \forall i \in \mathcal{N}
            \label{eq:DSDPOPF:qg}\\
        & 
            -\lambdaP_{i} - \lambdaPf_{ij} + \nuThermalPfr_{ij} = 0
            && \forall ij \in \mathcal{E}
            \label{eq:DSDPOPF:pf}\\
        & 
            -\lambdaQ_{i} - \lambdaQf_{ij} + \nuThermalQfr_{ij} = 0
            && \forall ij \in \mathcal{E}
            \label{eq:DSDPOPF:qf}\\
        & 
            -\lambdaP_{j} - \lambdaPt_{ij }+ \nuThermalPto_{ij} = 0
            && \forall ij \in \mathcal{E}
            \label{eq:DSDPOPF:pt}\\
        & 
            -\lambdaQ_{j} - \lambdaQt_{ij} + \nuThermalQto_{ij} = 0
            && \forall ij \in \mathcal{E}
            \label{eq:DSDPOPF:qt}\\
        &
            \mathcal{A}_{R}(\lambda, \mu, \nu) + \im \mathcal{A}_{I}(\lambda, \mu, \nu) + \Ssdp = 0
            &&
            \label{eq:DSDPOPF:W}
            \\
        & 
            \muPgMin, \muPgMax, \muQgMin, \muQgMax, \muWmMin, \muWmMax \geq 0
            && 
            \label{eq:DSDPOPF:non_negative} \\
        & 
            \nuThermalfr_{ij} = (\nuThermalSfr_{ij}, \nuThermalPfr_{ij}, \nuThermalQfr_{ij}) \in \mathcal{Q}^{3},
            \nuThermalto_{ij}  = (\nuThermalSto_{ij}, \nuThermalPto_{ij}, \nuThermalQto_{ij})\in \mathcal{Q}^{3}
            && \forall ij \in \mathcal{E}
            \label{eq:DSDPOPF:cone:nu} \\
        &
            \Ssdp \in \hermat^{|\mathcal{N}|}_{+}
            \label{eq:DSDPOPF:psd}
            && 
    \end{align}
\end{subequations}
\end{model}

%% file: tex/learning.tex
\section{The Dual SDP Proxy Architecture}
\label{sec:architecture}

The core contribution of the paper is a dual-feasible architecture for DSDP-OPF,
which leverages a new dual completion layer that specifically handles the positive semidefinite constraint \eqref{eq:DSDPOPF:psd}.

Algorithm \ref{algo:DCP} presents the proposed DCP architecture.
First (step 1), an initial model, e.g., a deep neural network, predicts dual variables $\lambdaP, \lambdaQ$, $\nuThermalPfr, \nuThermalQfr, \nuThermalPto, \nuThermalQto$.
Dual variables $\lambdaPf$, $\lambdaQf$, $\lambdaPt$, $\lambdaQt$ are then obtained (step 2) by substituting these predicted values into dual equality constraints \eqref{eq:DSDPOPF:pf}-\eqref{eq:DSDPOPF:qt}.
Steps 3 and 4 use the same argument as \cite{qiu2024dual} to recover variables $\nuThermalSfr$, $\nuThermalSto$ and $\muPgMin$, $\muPgMax$, $\muQgMin$, $\muQgMax$ in closed-form, using constraints \eqref{eq:DSDPOPF:cone:nu} and \eqref{eq:DSDPOPF:pg}-\eqref{eq:DSDPOPF:qg} combined with dual optimality conditions.
Finally, $\mathbf{S}$ and $\muWmMin, \muWmMax$ are recovered as follows.
Letting $\muWmMin = \muWmMax = 0$ in Eq. \eqref{eq:AR_def}, define
\begin{align}
    \label{eq:completion:S_hat}
    \hat{\mathbf{S}} = -\mathcal{A}_{R}(\lambda, \mu, \nu) - \im \mathcal{A}_{I}(\lambda, \mu, \nu),
\end{align}
and let $\mathbf{S} = \hat{\mathbf{S}} - \delta I$, where $\delta = \min(0, \lambda_{min}(\hat{\mathbf{S}})) \leq 0$.
Note that $\delta$ is well defined because $\hat{\mathbf{S}}$ is Hermitian by construction, and that $\mathbf{S} \succeq 0$.
Next, let $(\muWmMin_{i}, \muWmMax_{i}) = (0, -\delta), \forall i \in \mathcal{N}$, which immediately ensures that the completed dual solution satisfies constraint \eqref{eq:DSDPOPF:W}.

The dual completion strategy presented in Algorithm \ref{algo:DCP} ensures that the proposed DCP architecture always outputs dual feasible solutions.
Using a neural network as the prediction model, training can be done using gradients of the objective function \eqref{eq:DSDPOPF:obj} propagated through the completion layer in an end-to-end fashion.
Training the model in this self-supervised way by using the objective function as training loss removes the need to generate ground truth solutions, which can be time-consuming.
For a more detailed description of the DCP architecture, the reader may refer to \cite{qiu2024dual, tanneau2024dual}.

\begin{algorithm}[!t]
    \caption{The DCP Methodology for DSDP-OPF.}
    \label{algo:DCP}
    \setcounter{AlgoLine}{0}
    Predict $\lambdaP, \lambdaQ$, $\nuThermalPfr, \nuThermalQfr, \nuThermalPto, \nuThermalQto$ \\
    Recover $\lambdaPf, \lambdaQf, \lambdaPt, \lambdaQt$ using constraints \eqref{eq:DSDPOPF:pf}--\eqref{eq:DSDPOPF:qt} \label{algo:recover_nu} \\
    Recover $\nuThermalSfr, \nuThermalSto$ using constraint \eqref{eq:DSDPOPF:cone:nu} \\
    Recover $\muPgMin, \muPgMax, \muQgMin, \muQgMax$ via \eqref{eq:DSDPOPF:pg}, \eqref{eq:DSDPOPF:qg} \label{algo:recover_mu} \\
    Recover $\mathbf{S} = \hat{\mathbf{S}} - \delta I$, $\muWmMin=0$, $\muWmMax=-\delta$, where $\hat{\mathbf{S}}$ is obtained as per Eq. \eqref{eq:completion:S_hat} and  $\delta = \min(0, \lambda_{min}(\hat{\mathbf{S}}))$ \label{algo:recover_S}
\end{algorithm}

%% file: tex/results.tex
\section{Numerical Experiments}
\label{sec:results}

\newcommand{\ieeeXXXS}{\texttt{ieee14}}
\newcommand{\ieeeXXS}{\texttt{ieee30}}
\newcommand{\pegaseXS}{\texttt{pegase89}}
\newcommand{\ieeeXS}{\texttt{ieee118}}
\newcommand{\ieeeS}{\texttt{ieee300}}

\subsection{Experiment Details}
\label{sec:results:data_generation}

The paper conducts experiments on several systems from the PGLib library v21.07 \cite{babaeinejadsarookolaee2019power} with up to 500 buses.
AC-OPF, SOC-OPF and SDP-OPF instances are formulated using PowerModels \cite{power_models}, and solved with Ipopt \cite{wachter2006implementation} and Mosek 10 \cite{mosek}, respectively.
Both solvers use tolerances of $10^{-6}$, with other parameters set to default.
For each system, 20,000 OPF instances are generated by perturbing reference $\Sd$ values using the method in \cite{anonymous2024pglearn}.
Each dataset is then split into training (90\%), validation (5\%) and testing (5\%).
The solutions of the solvers are only used for evaluating optimality gaps during testing, and are not used in training nor validation.

The proposed SDP proxy is evaluated against the SOC proxy proposed in \cite{qiu2024dual}, which is based on the SOC relaxation of AC-OPF.
Both proxies use the same training, validation and test data.
The DCP proxies are implemented in Python 3.10 using PyTorch 2.0 \cite{paszke2019pytorch} and trained with the Adam optimizer \cite{kingma2014adam}.
All experiments are conducted on Intel Xeon 2.7GHz CPU machines running Linux and equipped with Tesla RTX 6000 GPUs.

\subsection{Performance Evaluation}
The quality of the DCP dual solutions is evaluated in terms of optimality gap with respect to the objective values of Mosek solutions (to the dual problems) and those of Ipopt solutions (to the original AC-OPF).

The quality of a dual bound is evaluated through its dual optimality gap, defined as 
\begin{equation}
    \text{dual optimality gap} = \frac{z^{*}_{AC} - z}{z^{*}_{AC}} \geq 0,
\end{equation}
where $z$ is a valid dual bound, obtained from a dual-feasible solution of the SOC or SDP relaxation, and $z^{*}_{AC}$ is the objective value of the best-known AC-OPF solution, obtained by Ipopt.
The paper also reports the gap closed by the proposed SDP proxy compared to the SOC proxy of \cite{qiu2024dual}, defined as
\begin{align}
    \label{eq:gap_closed}
    \text{gap closed} = \frac{\hat{z}_{SDP} - \hat{z}_{SOC}}{z^{*}_{SDP} - \hat{z}_{SOC}},  
\end{align}
where $z^{*}_{SDP}$, $\hat{z}_{SDP}$ and $\hat{z}_{SOC}$ denote the dual bound obtained from an optimal dual SDP solution, the predicted dual SDP solution, and the predicted dual SOC solution, respectively.
A positive (resp. negative) gap closed indicates that the SDP proxy yields better (resp. worse) dual bounds that the SOC proxy; a 100\% gap closed indicates the SDP proxy produces dual-optimal solutions.
The paper uses the geometric mean $\mu(x_{1}, ..., x_{n}) \, {=} \sqrt[n]{x_{1} x_{2} ... x_{n}}$ to report mean dual gaps, and the arithmetic mean to report mean gap closed, since the latter may take negative values.

\subsection{Dual Conic Proxy Performance}
\label{sec:results:DCP}

\begin{table}[!t]
    \centering
    \caption{DCP Performance Results}
    \label{tab:results:gaps}
    \begin{tabular}{lrrrrr}
        \toprule
        & \multicolumn{4}{c}{Dual Optimality Gap w.r.t. $z^{*}_{AC}$}\\
        \cmidrule(lr){2-5}
        System
            & \multicolumn{1}{c}{$z^{*}_{SOC}$}
            & \multicolumn{1}{c}{$z^{*}_{SDP}$}
            & \multicolumn{1}{c}{$\hat{z}_{SOC}$}
            & \multicolumn{1}{c}{$\hat{z}_{SDP}$}
            & $^{\dagger}$Gap closed
            \\
\midrule
\multirow{1}{*}{\ieeeXXXS}
        & 0.101 (0.009)
        & 0.000 (0.001)
        & \textbf{0.153} (0.050)
        & 0.404 (0.062)
        & -170.232 (54.157)
        \\
\multirow{1}{*}{\texttt{ieee30}}
        & 21.062 (2.490)
        & 0.029 (0.030)
        & 24.338 (3.074)
        & \textbf{23.758} (2.656)
        & 1.456 (12.851)
        \\
\multirow{1}{*}{\ieeeXS}
        & 0.837 (0.683)
        & 0.037 (0.110)
        & 2.110 (1.008)
        & \textbf{0.713} (0.700)
        & 61.261 (21.094)
        \\
\multirow{1}{*}{\ieeeS}
        & 2.009 (0.663)
        & 0.373 (0.915)
        & \textbf{5.580} (1.516)
        & 5.693 (1.687)
        & -2.601 (\phantom{0}4.761)
        \\
\multirow{1}{*}{\texttt{goc500}}
        & 0.097 (0.024)
        & 0.000 (0.000)
        & 0.384 (0.564)
        & \textbf{0.229} (0.582)
        & 34.326 (25.455)
        \\
        \bottomrule
    \end{tabular} \\
    \footnotesize{All values are in \%. Mean and standard deviation (in parentheses) are presented. $^{\dagger}$Gap closed by dual SDP proxy ($\hat{z}_{SDP}$) compared to SOC proxy ($\hat{z}_{SOC}$) with respect to ground truth SDP bound ($z^{*}_{SDP}$), see Eq \eqref{eq:gap_closed}.}
\end{table}

Table \ref{tab:results:gaps} presents the performance of the SOC and SDP proxies.
The proposed SDP proxy yields smaller dual gaps than the SOC proxy on the 30-, 118- and 500-bus systems.
It is important to note that the performance of the SOC proxy is intrinsically limited by the quality of the SOC relaxation, which can be significantly worse than the SDP relaxation, e.g., on the {\ieeeXXS} system ($>20\%$ vs $0.029\%$ gap).
In particular, on the 118-bus system,
\emph{the SDP proxy outperforms the SOC relaxation} which, in turn, implies that no SOC proxy can match the SDP proxy on this test case.

Furthermore, the results illustrate the challenges of training high-quality SDP proxy, especially on the {\ieeeXXS} system.
This is attributed to the higher output dimension of the SDP proxy, and the additional difficulty of the PSD constraint \eqref{eq:DSDPOPF:psd}.
The numerical results of Table \ref{tab:results:gaps} thus confirm the benefits of using the SDP relaxation as the basis for building high-quality dual proxies for AC-OPF.

Finally, timing results, not reported for lack of space, further demonstrate that the SDP proxy achieves several orders of magnitude speedups compared to Mosek.

%% file: tex/conclusion.tex
\section{Conclusion}
\label{sec:conclusion}

This work has proposed, for the first time, a dual-feasible proxy for the SDP relaxation of AC-OPF.
The proxy leverages a novel dual completion mechanism for SDP-OPF, is trained in a self-supervised fashion, and provides high-quality lower bounds for small- to medium-size systems significantly faster than conventional solvers.
Compared to the original proposed DCP that uses the SOC relaxation, it demonstrates notable improvements in optimality gap for some small- and medium-size systems.

Preliminary results on systems with up to 500 buses validate the potential of this approach, and encourage further investigation into applying the DCP methodology to SDP problems.
A promising avenue consists in scaling SDP-DCP to larger systems by exploiting chordal sparsity, which is key to making SDP optimization tractable in practice.